\documentstyle{article}
\begin{document}
\title{{Nonclassical   Solutions of Fully Nonlinear Elliptic
Equations II: Hessian  Equations  and Octonions  } }
 \author{{Nikolai Nadirashvili\thanks{LATP, CMI, 39, rue F. Joliot-Curie, 13453
Marseille  FRANCE, nicolas@cmi.univ-mrs.fr},\hskip .4 cm Serge
Vl\u adu\c t\thanks{IML, Luminy, case 907, 13288 Marseille Cedex
FRANCE, vladut@iml.univ-mrs.fr} }}

\date{}
\maketitle

\def\n{\hfill\break} \def\al{\alpha} \def\be{\beta} \def\ga{\gamma} \def\Ga{\Gamma}
\def\om{\omega} \def\Om{\Omega} \def\ka{\kappa} \def\lm{\lambda} \def\Lm{\Lambda}
\def\dl{\delta} \def\Dl{\Delta} \def\vph{\varphi} \def\vep{\varepsilon} \def\th{\theta}
\def\Th{\Theta} \def\vth{\vartheta} \def\sg{\sigma} \def\Sg{\Sigma}
\def\bendproof{$\hfill \blacksquare$} \def\wendproof{$\hfill \square$}
\def\holim{\mathop{\rm holim}} \def\span{{\rm span}} \def\mod{{\rm mod}}
\def\rank{{\rm rank}} \def\bsl{{\backslash}}
\def\il{\int\limits} \def\pt{{\partial}} \def\lra{{\longrightarrow}}

\section{Introduction}
\bigskip

This paper is a sequel to [NV1]; we study here a class of fully
nonlinear second-order elliptic equations of the form
$$F(D^2u)=0\leqno(1.1)$$
defined in a domain of ${\bf R}^n$. Here $D^2u$ denotes the
Hessian of the function $u$. We assume that
$F$ is a Lipschitz  function defined on the space $ S^2({\bf R}^n)$ 
of ${n\times n}$ symmetric matrices  satisfying the uniform ellipticity condition,
 i.e. there exists a constant $C=C(F)\ge 1$ (called an {\it ellipticity
constant\/}) such that 
$$C^{-1}||N||\le F(M+N)-F(M) \le C||N||\;
\leqno(1.2)$$ 
for any non-negative definite symmetric matrix $N$; if
$F\in C^1(S^2({\bf R}^n))$ then this condition is equivalent to
$${1\over C'}|\xi|^2\le F_{u_{ij}}\xi_i\xi_j\le C' |\xi |^2\;,
\forall\xi\in {\bf R}^n\;.\leqno(1.2')$$
 Here, $u_{ij}$ denotes the partial derivative
$\pt^2 u/\pt x_i\pt x_j$. A function $u$ is called a {\it
classical\/} solution of (1) if $u\in C^2(\Om)$ and $u$ satisfies
(1).  Actually, any classical solution of (1) is a smooth
($C^{\alpha +3}$) solution, provided that $F$ is a smooth
$(C^\alpha )$ function of its arguments.

For a matrix $S \in   S^2({\bf R}^n)$  we denote by $\lambda(S)=\{
\lambda_i : \lambda_1\leq...\leq\lambda_n\}
 \in {\bf R}^n$  the (ordered) set  of eigenvalues of the matrix $S$.
Equation (1) is called a Hessian equation  ([T1],[T2] cf. [CNS])
 if the function $F(S)$ depends only
  on the eigenvalues $\lambda(S)$ of the matrix $S$, i.e., if
 $$F(S)=f(\lambda(S)),$$
 for some function $f$  on ${\bf R}^n$ invariant under  permutations of
 the coordinates.

 In other words the equation (1) is called Hessian if it is invariant under
 the action of the group
 $O(n)$ on $S^2({\bf R}^n)$:
 $$\forall O\in O(n),\; F({^t O}\cdot S\cdot O)=F(S) \;.\leqno(1.3) $$
 The Hessian invariance relation (3) implies the following:

\medskip
 (a) $F$ is a smooth (real-analytic) function of its arguments if and only if $f$ is
a smooth (real-analytic) function.

\medskip
 (b) Inequalities (1.2) are equivalent to the inequalities
 $${\mu\over C_0} \leq { f ( \lambda_i+\mu)-f ( \lambda_i) } \leq C_0 \mu,
 \; \forall  \mu\ge 0,$$
 $\forall  i=1,...,n$, for some positive constant $C_0$.

\medskip
 (c) $F$ is a concave function if and only if $f$ is concave.

\medskip
 Well known examples of the Hessian equations are Laplace, Monge-Amp\`ere,
Bellman, Isaacs and Special Lagrangian equations.

 \medskip
 Bellman and Isaacs equations appear in the theory of controlled diffusion processes, see [F]. The
 both are fully nonlinear uniformly elliptic equations of the form (1.1). The Bellman equation
 is concave in $D^2u \in  S^2({\bf R}^n)$ variables. However, Isaacs operators are, in  general,
 neither concave nor convex. In a simple homogeneous form the Isaacs equation can be written
 as follows:
$$F(D^2u)=\sup_b \inf_a L_{ab}u =0, \leqno (1.4) $$
 where $L_{ab}$ is a family of linear uniformly elliptic operators of type 
 $$L=  \sum a_{ij} {\partial^2 \over  \partial x_i \partial x_j } \leqno (1.5)$$
  with an ellipticity
 constant $C>0$ which depends on two parameters $a,b$.

\medskip
Consider the Dirichlet problem
$$\cases{F(D^2u)=0 &in $\Om$\cr
u=\vph &on $\pt\Om\;,$\cr}\leqno(1.6)$$ where  $\Omega \subset {\bf
R}^n$ is a bounded domain with smooth boundary $\partial \Omega$
and $\vph$ is a continuous function on $\pt\Om$.

We are interested in the problem of existence and regularity of
solutions to  Dirichlet problem (1.6) for Hessian equations and Isaacs equation.
The problem (1.6) has always a unique viscosity (weak)
solution for  fully nonlinear elliptic equations (not necessarily
Hessian equations). The viscosity solutions  satisfy the equation
(1.1) in a weak sense, and the best known interior regularity
([C],[CC],[T3]) for them is $C^{1+\epsilon }$ for some $\epsilon
> 0$. For more details see [CC], [CIL]. Until recently it remained
unclear whether non-smooth viscosity solutions exist. In [NV1] we
proved   the existence in 12 dimensions  of non-classical
viscosity solutions to a fully nonlinear elliptic equation. The
paper [NV1] uses the function $$w_{12}(x)={Re (q_1q_2q_3)\over |x|},
$$ where $q_i\in {\bf H},\ i=1,2,3,$ are Hamiltonian quaternions,
$x=(q_1,q_2,q_3)\in {\bf H}^3={\bf R}^{12}$  which is a viscosity solution in
${\bf R}^{12}$ of a uniformly elliptic equation (1.1) with a smooth
$F$. Moreover, in [NV2] we proved that in 24 dimensions there exists 
a singular viscosity  solution to  a uniformly elliptic equation (1.1) with a smooth
$F$ which lies in $C^{2-\varepsilon}$ for a small positive $\varepsilon$.
 
\medskip
Our first main goal  is to show that  an octonionic analogue of
$w_{12}$  provides singular solutions  to Hessian
  uniformly elliptic equations in 21 (and more)
dimensions. Moreover the following theorem holds for a certain 
harmonic cubic polynomial $P_{24}$ in ${\bf R}^{24}$:

\bigskip

{\bf Theorem 1.1.}

{\it For any $\delta , \; 1\leq \delta < 2$  and  any plane
$H'\subset{\bf R}^{24}, \; \dim H'=21$ the function
$$(P_{24} (x)/ |x|^{\delta })_{|H'}
$$ is a viscosity solution  to  a uniformly elliptic Hessian equation
  $(1.1)$ in a unit ball $B\subset {\bf
R}^{21}$  for  the cubic form
$$ P_{24}(x)={Re((o_1\cdot o_2)\cdot o_3)}={Re(o_1\cdot(o_2\cdot o_3))},
$$ where $o_i\in  {\mathcal{O}},\; i=1,2,3,$  $\mathcal{O}$ being
the  algebra of Caley octonions, $x=(o_1,o_2,o_3)\in {\mathcal{O}}^3={\bf
R}^{24}.$}

\medskip
\noindent It shows the optimality of the result by Caffarelli-Trudinger [C,CC,T3] on
 the interior $C^{1,\alpha}$-regularity of viscosity solutions of fully nonlinear 
equations,  even in the Hessian case.

\medskip
The second main goal  is to show that the same function  is a viscosity solution  to
a uniformly elliptic Isaacs equation:
\bigskip

{\bf Theorem 1.2.}

{\it For any $\delta , \; 1\leq \delta < 2$  and  any plane
$H'\subset{\bf R}^{24}, \; \dim H'=21$ the function
$$(P_{24} (x)/ |x|^{\delta })_{|H'}
$$ is a viscosity solution  to  a uniformly elliptic  Isaacs equation
  $(1.4)$ in a unit ball $B\subset {\bf
R}^{21}$.}

\medskip
 The rest of the paper is organized as follows: in Section 2 we recall some preliminary results, 
 we introduce the form $P_{24}$ and give its main properties in Section 3, we prove 
Theorem 1.1 in Section 4,  and, finally, we prove Theorem 1.2 in Section 5.

\section{Preliminary results }

\medskip
   Let $w=w_\delta$ be a homogeneous function
of order $3-\delta ,\  1\leq \delta <2 $, defined on a unit ball
$B =B_1\subset {\bf R}^n$ and smooth in $B \setminus\{0\}$. Then
the Hessian of $w$ is homogeneous of order $(1-\delta)$. Define
the map
$$\Lambda :B  \longrightarrow \lambda (S) \in {\bf R}^n\; .$$
    $\lambda(S)=\{ \lambda_i :
\lambda_1\leq...\leq\lambda_n\}
 \in {\bf R}^n$ being  the (ordered) set  of eigenvalues of the matrix
 $S=D^2w$.

Let $K\subset {\bf R}^n$ be an open convex cone, such that
$$ \{ x\in {\bf R}^n : x_i\geq 0,\ i=1,...,
n\}  \subset K .$$

Set
$$L:={\bf R}^n\setminus (K\cup -K).$$

  We say that a set  $E\subset {\bf R}^n$ satisfy $K$-cone condition if
  $\;  (a-b)\in L$
for any $ a,b\in E.$

Let $\Sigma_n$ be the group  of permutations of   $\{ 1,...,n\}$. For
any $\sigma \in \Sigma_n$, we denote by $T_{\sigma}$ the linear
transformation of ${\bf R}^n$ given by $x_i \mapsto x_{\sigma(i)}, \; i=1,...,n.$

\bigskip\noindent {\bf  Lemma 2.1.} {\it  Assume that
$$M:=\bigcup_{\sigma \in \Sigma_{n} }\   T_{\sigma }\Lambda (B)\subset {\bf R}^n $$
satisfies the $K$-cone condition. If $\delta> 0$ we assume additionally that $w$
 changes sign in $B$. Then $w$ is a viscosity solution in $B$ of a uniformly 
elliptic Hessian equation $(1)$.}

\bigskip
\noindent
{\em   Proof  }. Let us choose in the space ${\bf R}^n$ an
orthogonal coordinate system $z_1,\dots,z_{n-1},s,$  such that
$s=x_1+...+x_n$ . Let $\pi : {\bf R}^n\to Z$ be the orthogonal
projection of ${\bf R}^n$ onto the $z$-space. Let $K^\ast$ denote
the adjoint cone of $K$, that is, $K^\ast = \{b\in {\bf R}^n:
b\cdot c \ge 0 \ for\ all\ c \in K\} $. Notice that $a \in L $
implies $a\cdot b =0$ for some $b \in K^\ast$.  We represent the
boundary of the cone $K$ as the graph of a Lipschitz function
$s=e(z)$, with $e(0)=0$, function $e$ is smooth outside the
origin: $$e(z)\ =\inf\{ c:\ (z+cs)\in K \} .$$

Set $m=\pi \bigl (M )$. We prove that $M$ is a graph of a
Lipschitz  function on $m$, $$M =\{z\in m:s=g(z)\}\;.$$
Let $a,\hat a \in M ,  a =(z,s),\hat a = (\hat z,\hat s)$. Since
$a-\hat a \in L$, we have $$-e(z-\hat z) \le \hat s - s \le e(z-
\hat z).$$ Since $e(0)=0, g(z):=s$ is single-valued. Also
$$|g(z)-g(\hat z)|=|s-\hat s|\le |e(z-\hat z)|\le C|z-\hat z|.$$

The function $g$ has an extension $\widetilde{g}$ from the set $m$
to ${\bf R}^{n-1}$ such that $\widetilde{g}$ is a Lipschitz
function and the graph of $\widetilde{g}$ satisfies the $K$-cone
condition. One can define such extension $\widetilde{g}$ simply by
the formula
$$\widetilde{g}(z):=\inf_{w\in m } \bigl\{ g(w)\ +  e(z-w)
\bigr\}\;.$$

To show that this formula works let $(z, \tilde g(z)),  (\hat
z,\tilde g(\hat z))$ lie in the graph $\tilde g$. We must show
$$-e(z-\hat z)\ \le \ \tilde g(z)-\tilde g(\hat z)\ \le e(z-\hat z).$$
Now
$$\tilde g(\hat z)\ =\ g(w)+e(\hat z-w)$$
for some $w\in m$. Thus
$$\tilde g(z)-\tilde g(\hat z)\ \le \ g(w)+e(z-w)-(g(w)+e(\hat z-w))\ \le \
e(z-\hat z),$$ since $e(a+b)\ \le \ e(a)+e(b)$, as $e(\cdot )$ is
convex, homogenous. Similarly
$$\tilde g(z)-\tilde g(\hat z)  \ge  -e(z-\hat z).$$
\medskip
Let us set
$$
f':= s - \widetilde{g}(z).
$$
Since the level surface of the function  $ f'$ satisfies $K$-cone
condition it follows that $\nabla f \in K^*$ a. e.  where  $K^*$
is the adjoint cone to $K$. Moreover   the function $w$ satisfies
the equation
$$
 f'(\lambda (S))= 0.
$$
on $B \setminus \{0\}$.

Set
$$f=\sum_{\sigma \in {\Sigma_n}} f'(\sigma (x)).$$
Then $f$ is a Lipschitz  function invariant under the action of
the group $\Sigma_n$ and satisfies the equation
$$
 f(\lambda (S))= 0.
$$
on $B \setminus \{0\}$.

\medskip
We show now that $w$ is a viscosity solution of (1) on the whole
ball $B$.

Assume first that $\delta =1$. Let $p(x),\ x\in B$ be a quadratic
form such that $p\leq w$ on $B$. We choose any quadratic form
$p'(x)$ such that $p\leq p'\leq w$ and there is a point $x'\neq 0$
at which $p'(x')= w(x')$. Then it follows that $F(p)\leq F(p')\leq
0$. Consequently for any quadratic form $p(x)$ from the inequality
$p\leq w$ ($p\geq w$) it follows that $F(p)\leq 0$ ($F(p)\geq 0$).
This implies that $w$ is a viscosity solution of (1) in $B$ (see
Proposition 2.4 in [CC]).

If $1< \delta < 2$ then for any smooth function $p$ in  $B$ the
function $w-p$ changes sign in any neighborhood of $0$. Hence, by
the same proposition in [CC], it follows that $w$ is a viscosity
solution of (1) in $B$.
\bigskip

Next we need the following property of the eigenvalues $
\lambda_1\ge\lambda_2\ge\ldots\ge\lambda_{n}$ of real symmetric matrices 
 of order $n$ which is a classical result by Hermann  Weyl [We]:

\bigskip\noindent {\bf  Lemma 2.2.} 
  {\em Let $ A,B$  be two real symmetric  matrices
with the eigenvalues $
\lambda_1\ge\lambda_2\ge\ldots\ge\lambda_{n} $ and $
\lambda'_1\ge\lambda'_2\ge\ldots\ge\lambda'_{n} $ respectively.
Then for the eigenvalues $
\Lambda_1\ge\Lambda_2\ge\ldots\ge\Lambda_{n} $ of the matrix $A-B$
we have}
$$ \Lambda_1\ge\max_{i=1,\cdots, n}(\lambda_i-\lambda'_i),
 \;\;\Lambda_n\le\min_{i=1,\cdots, n}(\lambda_i-\lambda'_i).$$

\medskip
We need also the following simple fact:

\medskip
{\bf Lemma 2.3. }
{\em Let $L:{\bf R}^{n}\longrightarrow {\bf R}^{n} $ be a
symmetric linear operator with  the eigenvalues $\lambda_1\ge
\lambda_2\ge\ldots\ge\lambda_{n} $ and let  $H$ be a hyperplane
$H\subset {\bf R}^{n}$  invariant under $L.$ Then for the
eigenvalues $\lambda_1'\ge \lambda_2'\ge\ldots\ge\lambda_{n-1}' $
of the restriction   $L_{|H}$ one  has}
$$\lambda_1\ge\lambda_1'\ge
\lambda_2\ge \lambda_2'\ge\ldots\ge\lambda_{n-1}\ge
\lambda_{n-1}'\ge\lambda_{n}.
$$

\section{Cubic form $P=P_{24} $ }
In this section we introduce and investigate the cubic form which
will be used to construct our singular solutions. It is based on 
 the  algebra of Caley octonions $\mathcal{O}$; for this  algebra
 we use the notation and conventions in [Ba] (in particular, $e_1 e_2 = e_4$).
 Let $ V=(X,Y,Z)\in {\bf R}^{24}$ be a variable vector with $X,Y,$ and $Z\in {\bf R}^8.$
For any $ t=(t_0,t_1,\ldots,t_8)\in {\bf R}^8$ we denote  by
$$
ot=t_0+t_1\cdot e_1+t_2\cdot e_2+\ldots+t_7\cdot e_7
 \in \mathcal{O}$$
 its natural image in $\mathcal{O}.$  For any
  $ o =o_0+o_1\cdot e_1+o_2\cdot e_2+\ldots+o_7\cdot e_7\in \mathcal{O}$
 its conjugate
will be denoted $ o^*=o_0-o_1\cdot e_1-o_2\cdot
e_2-\ldots-o_7\cdot e_7$; thus, $o^*\cdot o=o\cdot o^*=\mid
o\mid^2.$

    Define the cubic form $P=P_{24}(V)=P(X,Y,Z) $ as follows
    $$ P(X,Y,Z)=Re((oX\cdot oY)\cdot oZ)=Re(oX\cdot (oY\cdot oZ))=$$ $$
    (Z_0Y_0-Z_1Y_1-Z_2Y_2-Z_3Y_3-Z_4Y_4-Z_5Y_5-Z_6Y_6-Z_7Y_7)X_0+$$ $$
     (-Z_1Y_0-Z_0Y_1-Z_4Y_2-Z_7Y_3+Z_2Y_4-Z_6Y_5+Z_5Y_6+Z_3Y_7)X_1+$$
     $$ (-Z_2Y_0+Z_4Y_1-Z_0Y_2-Z_5Y_3-Z_1Y_4+Z_3Y_5-Z_7Y_6+Z_6Y_7)X_2+$$
     $$ (-Z_3Y_0+Z_7Y_1+Z_5Y_2-Z_0Y_3-Z_6Y_4-Z_2Y_5+Z_4Y_6-Z_1Y_7)X_3 +$$
     $$(-Z_4Y_0-Z_2Y_1+Z_1Y_2+Z_6Y_3-Z_0Y_4-Z_7Y_5-Z_3Y_6+Z_5Y_7)X_4 +$$
     $$(-Z_5Y_0+Z_6Y_1-Z_3Y_2+Z_2Y_3+Z_7Y_4-Z_0Y_5-Z_1Y_6-Z_4Y_7)X_5 + $$
     $$(-Z_6Y_0-Z_5Y_1+Z_7Y_2-Z_4Y_3+Z_3Y_4+Z_1Y_5-Z_0Y_6-Z_2Y_7)X_6+$$
      $$(-Z_7Y_0-Z_3Y_1-Z_6Y_2+Z_1Y_3-Z_5Y_4+Z_4Y_5+Z_2Y_6-Z_0Y_7)X_7.$$

\medskip
Its principal property for us is

\medskip
     {\bf Proposition 3.1. } {\em  Let $a=(x,y,z)\in S_1^{23};$   define
$$
W=W(a)=P(a),\; m=m(a)=m(x,y,z)= |x|\cdot |y| \cdot |z| .$$

 Then the characteristic polynomial $CH(T)=CH_{P,a}(T)$ of the
Hessian $H(a)= D^2P(a)$ is given by $$CH(T)\;
=(T^3-T+2m)(T^3-T-2m)(T^3-T+2W)^6. $$}
\medskip

\vskip .3 cm {\em Proof.} The weak associativity  $ Re((oX\cdot oY)\cdot oZ)=Re(oX\cdot (oY\cdot oZ)) $ is Corollary 15.12, p.110 of the book [Ad]. Proposition 5.7  [Ad, p.35] and Theorem 15.14 [Ad, p.111] show that the triality polynom $P(X,Y,Z)$ is Spin(8)-invariant. Thus  the characteristic polynomial $CH(T)$ is invariant under the action of Spin(8), and we can suppose (applying  the action) that the vectors $x\in {\bf R},\; y\in {\bf R}+e_1{\bf R}, \; z\in {\bf R}+e_1{\bf R}+e_2{\bf R}\subset {\mathcal{O}};$ thus $x,y,z\in {\bf H}\subset {\mathcal{O}}$ where ${\bf H}$ is generated by $ \{1,e_1,e_2 ,e_4\}.$  
 Brute force calculations give for the Hessian of $P$ relatively to the following ordering of
coordinates in ${\bf R}^{24}$:
$$\{X_0,X_1,X_2,X_4,Y_0,Y_1,Y_2,Y_4,Z_0,Z_1,Z_2,Z_4,
X_5,X_6,X_3,X_7,Y_5,Y_6,Y_3,Y_7,Z_5,Z_6,Z_3,Z_7\}$$

$$ H(a)=     \left(%
\begin{array}{cc}H_0
    & 0  \\
   0&H_1
\\
\end{array}%
\right) \quad\quad\quad\quad$$ for the following matrices
$H_0,H_1\in Mat_{12}({\bf R}):$

$$H_0=     \left(%
\begin{array}{ccc}0_4
    &  M_z & M_y \\
    ^t M_z&0_4 &M_x   \\
  ^tM_y & ^tM_x & 0_4
\\
\end{array}%
\right), \quad\quad  H_1=     \left(%
\begin{array}{ccc}0_4
    &  L_z & L_y \\
    ^t L_z&0_4 &L_x   \\
  ^tL_y & ^tL_x & 0_4
\\
\end{array}%
\right) \quad\quad\quad\quad $$ where
$$M_s=
  \left(%
\begin{array}{cccc}
  s_0&-s_1&-s_2&-s_3 \\
   -s_1&-s_0&-s_3&\; s_2 \\
  -s_2& s_3&-s_0&-s_1 \\
   -s_3&-s_2& s_1&-s_0\\
\end{array}%
\right), \quad L_s=
  \left(%
\begin{array}{cccc}
  -s_0&-s_1&s_2&-s_3 \\
   s_1&-s_0&-s_3&\; -s_2 \\
  -s_2& s_3&-s_0&-s_1 \\
   s_3&s_2& s_1&-s_0\\
\end{array}%
\right)$$  for an arbitrary $s=(s_0,s_1,s_2,s_3)\in {\bf R}^4$.

Direct  easy calculations show that $M_s ,L_s$ have the
following properties:

\medskip
1).$\;\; M_s\cdot {^tM}_s={^tM}_s \cdot M_s=  L_s\cdot
{^tL}_s={^tL}_s \cdot L_s= \mid s\mid^2 I_4;$

\noindent thus, $M_s ,L_s$  are proportional to   orthogonal
matrices. In particular, if  $\mid s\mid=1$    then $M_s ,L_s$ are
orthogonal themselves. We write $M_s= \mid s\mid O_s, \;L_s= \mid
s\mid O'_s$  with $O_s,O'_s\in O(4).$

\medskip
2).  $\det (M_s)=-\mid s\mid^4, \quad \det (O_s)=-1, \; \det
(L_s)=\mid s\mid^4, \quad \det (O'_s)=1;$

\medskip
3). the characteristic polynomials $ PM_s(T), PL_s(T)$ of $M_s,
L_s $ factor as
$$ PM_s(T)=(T^2-|s|^2)(T^2+2 s_0 T +|s|^2),\; PL_s(T)= (T^2+2 s_0 T + |s|^2)^2$$
  and those of $ O_s,
O'_s$   as
$$
PO_s(T)=(T^2-1)(T^2+2 s^*_0 T +1),\;PO'_s(T)=(T^2+2 s^*_0 T+1)^2$$
with $s^*_0=s_0/| s|; $

\medskip
4). define the symmetric matrices $N_s=  (O_s+ {^tO}_s),\;N'_s=
(O'_s+ {^tO}'_s)  ;$  then their spectrums are
     $$ Sp(N_s)=\{ 2,-2,-2 s^*_0,-2 s^*_0\},\;
     Sp(N'_s)=\{ -2s^*_0,-2s^*_0,-2 s^*_0,-2 s^*_0\};$$

     \medskip
5). For the product matrices  $ M_{rst}=M_r\cdot M_s\cdot M_t,\;
L_{rst}=L_r\cdot L_s\cdot L_t$, $r,s,t \in {\bf R}^4$ we have the
characteristic polynomials  $PM_{rst}, PL_{rst}$ of $
M_{rst},L_{rst}$:
$$  PM_{rst} (T)  =(T^2-| r|^2
 |s|^2| t|^2)
(T^2+2P(r,s,t)T+| r|^2 | s|^2|t|^2),\; $$
$$PL_{rst} (T)  = (T^2+2P(r,s,t)T+|r|^2 |s|^2|t|^2)^2.
$$

Let us calculate the characteristic polynomial $F$ of  $H_0$, the characteristic
 polynomial $G$  of $H_1$ being calculated in the same way using $L_s$ instead
  of $M_s$. Conjugating $H_0$ by the orthogonal
matrix $$\left(%
\begin{array}{ccc}^tO_z
    &  0_4 & 0_4 \\
    0_4&I_4 &  0_4 \\
  0_4 & 0_4& O_x
\\
\end{array}%
\right)$$ one gets

$$\tilde {H_0}=     \left(%
\begin{array}{ccc}0_4
    & |z| I_4 & \quad| y| ^tO_{xyz} \\
   \quad| z| I_4&0_4 & | x| I_4     \\
   \quad\quad| y| O_{xyz} & | x| I_4 & 0_4
\\
\end{array}%
\right) \quad\quad\quad\quad\quad\quad\quad\quad\quad\quad$$

\par Let now $\lambda\in Sp(\tilde {H_0})$, $v _\lambda=
(p_\lambda,q_\lambda,r_\lambda)$ being a corresponding
eigenvector, normalized by the condition  $ | v _\lambda |=1$.
 The condition $\tilde H_0\cdot
v_\lambda=\lambda v_\lambda$ gives
$$
\lambda p_\lambda=| z| q_\lambda+| y| {^tO}_{xyz}
r_\lambda\quad\quad\quad\quad\quad\quad\quad\quad
$$
$$
\lambda q_\lambda=| z| p_\lambda\quad+\quad| x|
r_\lambda\quad\quad\quad\quad\quad\quad\quad\quad
$$
$$
\lambda r_\lambda=| y| O_{xyz}  p_\lambda+| x|
q_\lambda\quad\quad\quad\quad\quad\quad\quad\quad.
$$
Multiplying the second and the third equations by $\lambda$ and
inserting in thus obtained equations the first one we get
$$
(\lambda^2-| z|^2)  p_\lambda=(| x|\cdot| z| + \lambda | y|
^tO_{xyz} )  r_\lambda \quad\quad\quad\quad\quad\quad\quad\quad
$$
$$
( \lambda^2-| x|^2)  r_\lambda=(| x|\cdot| z| + \lambda | y|
O_{xyz})  p_\lambda \quad\quad\quad\quad\quad\quad\quad\quad
  $$
which implies
  $$
  ( \lambda^2-| x|^2)(
\lambda^2-| z|^2) p_\lambda=(| x|\cdot| z| + \lambda | y|
^tO_{xyz})(| x|\cdot|z| + \lambda | y| O_{xyz}) p_\lambda
$$
and, after simplifying,
$$
\lambda (\lambda ^3I_4-\lambda I_4-mN_{xyz})p_\lambda=0 ,
$$
since $| x|^2+| y|^2+| z|^2=1$, $ m= | x|\cdot | y| \cdot | z|,$
$O_{xyz} {^tO}_{xyz}=I_4$,
$N_{xyz}=O_{xyz}+ {^tO}_{xyz}\; .$\\
  Hence, either $\lambda =0$ or
   $$
    (\lambda ^3 -\lambda)\in m\cdot Sp(N_{xyz})=\{-2m, 2m,-2W ,-2W \}.
$$
  This finishes the proof for $\lambda \neq 0$. If $\lambda =0$
  we get the conditions

$$
0=| z| q_\lambda+| y| ^tO_{xyz}
r_\lambda\quad\quad\quad\quad\quad\quad\quad\quad
$$
$$
0=| z| p_\lambda\quad+\quad| x|
r_\lambda\quad\quad\quad\quad\quad\quad\quad\quad
$$
$$
0=| y| O_{xyz}  p_\lambda+|x|
q_\lambda\quad\quad\quad\quad\quad\quad\quad\quad.
$$
immediately  implying  that $m=0$ (since else these equations give
$p_\lambda=0$) and the formula holds for this case as well.

\vskip .3 cm {\em Remark 3.1.}  If we do not instist on a computer-free
proof of the fact,  the inclusions   $x\in {\bf R},\; y\in {\bf R}+e_1{\bf R}, \; z\in {\bf R}+e_1{\bf R}+e_2{\bf R}$ will suffice. Indeed, the MAPLE instructions
 ($v$ being the coordinate vector)
$$H:=hessian(P,v):X2:=0:X4:=0:Y2:=0:Y4:=0: Z2:=0:Z4:=0:$$
$$X5:=0:X6:=0:X3:=0:X7:=0:Y5:=0:Y6:=0:Y3:=0: Y7:=0:$$            

\noindent $Z5:=0:Z6:=0:Z3:=0:Z7:=0:CH:= factor (charpoly(H,T));$
 \vskip .2 cm

\noindent return the formula of Proposition 3.1 in 20 seconds, $<60$ MB of space.

\vskip .3 cm The result of Proposition 3.1 can be written as
\vskip .3 cm {\bf Corollary 3.1.} {\em Define the angles
$\alpha,\beta $  by $ \alpha:=\arccos(3\sqrt 3 m),\;
\beta:=\arccos(3\sqrt 3 W) $. Then}
$$
Sp( H(a))  =\{ {2\over\sqrt 3} \cos(\alpha/3+\pi k/3),
6\times\{{2\over\sqrt 3}\cos(\beta/3+\pi(2l+1) /3)\}  \},
$$
 for $k=0,1,\ldots 5,\; l=0,1,2.$
  \vskip .3cm
  {\em Proof}.
  Indeed, if we put $\lambda={2\over\sqrt 3}\cos\gamma$, the equations
$\lambda ^3 -\lambda+2m=0$, $\lambda ^3 -\lambda-2m=0$ and
    $\lambda ^3 -\lambda+2W=0$
  become respectively,
  $\cos(3\gamma)=-\cos\alpha$,
$\cos(3\gamma)= \cos\alpha$ and $\cos(3\gamma)=-\cos\beta$ which
implies the result.

\medskip
  Let us  order the eigenvalues of $H(a)$ in the decreasing order:
$$\lambda_1\ge\lambda_2\ge\ldots \ge\lambda_{23}\ge\lambda_{24}.$$
 Since $|W|\le m $ and the cosine decreases on $[0,\pi] $ we get

\medskip {\bf Corollary 3.2.} {\em $$ \lambda_1={2\over\sqrt 3} \cos({\alpha\over 3}),
\lambda_2=\ldots=\lambda_{7}=\mu_1, \lambda_8=l_1 ,\lambda_9=l_2,
\lambda_{10}=\ldots=\lambda_{15}=\mu_2,$$
$$\lambda_{16}=-l_2,\lambda_{17}=-l_1,
\lambda_{18}=\ldots=\lambda_{23}=\mu_3, \lambda_{24}=-{2\over\sqrt
3} \cos({\alpha\over 3})$$ for

 $$l_1=\max\left\{{2\over\sqrt3}\cos({\alpha+\pi\over 3}  )
,{2\over\sqrt3}\cos({\alpha+5\pi\over 3})\right\},$$
 $$l_2=\min\left\{{2\over\sqrt3}\cos({\alpha+\pi\over 3})
,{2\over\sqrt3}\cos({\alpha+5\pi\over 3})\right\},$$
$\mu_1\ge\mu_2\ge\mu_3 $ being the roots of $T^3-T+2W=0$.}

 \medskip
 {\em Remark 3.2.}  We have the inequalities
 $$2\lambda_{3}\ge\lambda_1,\;\;\;2\lambda_{n-2}\le\lambda_{n}, n=12 \;\hbox{or} \; 24$$
which hold for the eigenvalues of $P_{24}$ as well for the form $P_{12}$ 
used in  [NV1]. They are essential for the  proofs in [NV1] and  are in fact 
the best possible. Indeed, one has the following result:

\medskip
 {\bf Proposition 3.2. } {\it Let $P\neq 0$ be a cubic form in ${\bf R}^n$.
 Then for some unit vector $d\in S^{n-1}_1\subset {\bf R}^n$
   the eigenvalues  $\lambda_1 \geq \lambda_2 \geq ...\ge\lambda_{n}$
   of the quadratic form $P_d:=\sum_i d_i P_{x_i}$ satisfy $$\lambda_1
\geq 2\lambda_2,\;\; 2\lambda_{n-1}\ge\lambda_{n}. $$   }

{\bf Proof }. Assume that at the point $a\in S^{n-1}_1 $ the cubic
form $P$ attains its supremum over $S^{n-1}_1$. Since $P$ is an odd
function on ${\bf R}^n$, $P(a)>0$. Choose  $d=a$ and let
$x_1,...,x_n$ be an orthonormal basis in ${\bf R}^n$ such that
$x_1$ is  directed along $d$. Since   the form $P$ attains at $d$
its supremum over $S^{n-1}_1$ it follows that in the coordinates
$x_i$ the cubic form $P$ contains no monoms  of the form
$cx_1^2x_i,\ i>1$. Thus the quadratic form $P_d$ contains no
monoms of the form $cx_1x_i,\ i>1$ and hence the vector $d$ is an
eigenvector of the quadratic form $P_d$ with the eigenvalue
denoted by  $\lambda $. Let $\lambda' $ be the maximal eigenvalue
of $P_d$ on the orthogonal complement of $d$ attained on the
eigenvector $b\in S^{n-1}_1$. The lemma will follows if we prove
that $\lambda \geq 2\lambda' $. We assume without loss that
$\lambda =1$ and that $x_2$ is directed along $b$. Then the
restriction of $P_d$ on the plane $\{ x_1,x_2 \} $ can be written
in the form
$$x_1^2+\lambda'x_2^2 $$
and  thus the restriction of the cubic form $P$ on this plane
becomes
$$x_1^3/3+\lambda'x_1x_2^2+cx_2^3.$$
It is easy to see that if $\lambda'>1/2$ then the supremum of the
function $P$ on the circle $x_1^2+x_2^2=1$ is not at the point
$(1,0)$ which implies the result.

\section{Proof of Theorem 1.1}
\medskip

Lemma 2.3 and Corollary 3.2 give

\medskip
{\bf Corollary 4.1.} {Let $a=(x,y,z)\in S_1^{23},$ let
$H=H_{19}\subset {\bf R}^{24}={\mathcal{O}}^3$  be a plane,
$dim(H)=19$ and let
$$\lambda'_1\ge\lambda'_2\ge\ldots
\ge\lambda'_{18}\ge\lambda'_{19}$$  be  the eigenvalues of the
Hessian  $D^2P_{|H}(a)$ written in the decreasing order. Then $$
  \lambda'_2
 =\mu_1,
\lambda'_{10} =\mu_2, \lambda'_{18}=\mu_3,
 $$
$\mu_1\ge\mu_2\ge\mu_3 $ being the roots of $T^3-T+2W=0$.}

\bigskip
 {\bf Proposition 4.1.} {\em Let $H\subset {\bf R}^{24},\;\dim H=21.$
 Set $M_\delta(u)={D^2w_{\delta}} _{|H}(u) $
for $u\in H,\; 1\leq\delta <2$.
 Suppose that  $a\neq b\in    H $  and
 let $O\in {\hbox {O}}({21} )$ be an orthogonal matrix s.t.
 $M_\delta(a,b,O):=M_\delta(a)- {^tO}\cdot M_\delta(b)\cdot O\neq 0$.
 Denote $ \Lambda_1\ge\Lambda_2\ge
 \ldots\ge\Lambda_{21}$  the eigenvalues of the matrix
 $M_\delta(a,b,O).$
  Then

$$ \varepsilon \le {\Lambda_1\over
-\Lambda_{21}}\le\varepsilon^{-1}
  $$
  for $\varepsilon:=\min\{ {2-\delta\over 4+\delta},{1\over 20}\}. $ }

\medskip
\noindent{\em Proof.}
We can suppose  without loss that $ |a| \le  |b|$, moreover, by
homogeneity we can suppose that $a\in S_1^{20}$ and thus  $ |b|
\ge 1.$ Let $\bar b:={b/|b|}\in S_1^{20}$ then
$M_{\delta}(b)=M_{\delta}(\bar b)|b|^{1-\delta} .$ One needs then
the following  result for the points $a,\bar b \in S_1^{20}:$

 \bigskip
{\bf Lemma 4.1. }{\em  Let  $\delta\in[1,2),$
   $a, \bar b\in S_1^{20},\; W=W(a),\;\bar W=W(\bar b),\; $ and let
$$\mu_1(\delta)=
 {2\over \sqrt
3}\cos\left({\arccos(3\sqrt 3 W)-\pi\over 3 }\right) -W\delta\ge$$
$$\mu_2(\delta)=
 {2\over \sqrt
3}\cos\left({\arccos(3\sqrt 3 W)+\pi\over 3 }\right) -W\delta\ge$$
$$
  \mu_3(\delta)=
 -{2\over \sqrt
3}\cos\left({\arccos(3\sqrt 3 W) \over 3}\right) -W\delta
  $$ $($resp.,
  $\bar\mu_1(\delta)\ge\bar\mu_2(\delta)\ge\bar\mu_3(\delta)\;)$
   be the roots of the polynomial
  $$P_{1,\delta}(T,W):=Q_1(T+\delta W)=$$
  $$T^3+3W\delta T^2+(3W^2\delta^2-1)T+W(2-\delta)+W^3\delta^3$$
$($resp. of the polynomial
$$ \bar P_{1,\delta}(T,\bar W):=Q_1(T+\delta\bar W)=$$
$$T^3+3\bar W \delta T^2+(3 \bar W^2\delta^2-1)T+\bar W(2-\delta)+\bar
 W^3\delta^3\; ).$$ Then  for any $K>0$ verifying
$ |K-1|+|\bar W-W|\neq 0$ one has
$${2-\delta\over 4+\delta}=:\varepsilon\le { \mu_+(K)\over -\mu_-(K)}\le
 {1\over\varepsilon}= {4+\delta\over 2-\delta}$$
where
$$\mu_-(K):=  \min\{\mu_1(\delta)-K\bar\mu_1(\delta),\;
\mu_2(\delta)-K\bar\mu_2(\delta),\;
\mu_3(\delta)-K\bar\mu_3(\delta)\} ,$$
$$\mu_+(K):=  \max\{\mu_1(\delta)-K\bar\mu_1(\delta),\;
\mu_2(\delta)-K\bar\mu_2(\delta),\;
\mu_3(\delta)-K\bar\mu_3(\delta)\}\;.  $$   }

\medskip{\em Proof of Lemma 4.1.} In the proof we will repeatedly use the
following elementary fact:
\medskip

{\em Claim. Let $ l_1\ge l_2 \ge l_3,$ $ l_1+ l_2+ l_3=t\ge 0,$  $
l_3\le -ht,$ with $h>0.$ Then $-l_1/l_3\in [h/(2h+1),(2h+1)/h] $
for $t>0$, $-l_1/l_3\in [1/2,2]$ for $t=0$.}
\medskip

 If $W= \bar W , K=1$ there is nothing to prove. If
$K=1$ one can suppose that $W>\bar W; $ we have
  $$(\mu_1(\delta)-K\bar\mu_1(\delta))+
(\mu_2(\delta)-K\bar\mu_2(\delta))+
(\mu_3(\delta)-K\bar\mu_3(\delta))=3(\bar W -W)\delta$$ and
$$ \mu_2(\delta)-K\bar\mu_2(\delta)=
 {2\over \sqrt 3}\left(\cos\left({\arccos(3\sqrt 3 W)+ \pi\over 3 }\right )-
 \cos\left({\arccos(3\sqrt 3 \bar W)+\pi\over 3 }\right)\right)$$
$$-(  W -\bar W)\delta\ge (2-\delta)(  W-\bar W).$$
  Therefore, one can take $\varepsilon=(2-\delta)/(4+\delta) $
  in this case.
 We can suppose then $W>\bar W, K\neq 1. $ Using the
relations $$\mu_1(\delta)(-W)=-\mu_3(\delta)(W),\;
\mu_2(\delta)(-W)=-\mu_2(\delta)(W),\;
\mu_3(\delta)(-W)=-\mu_1(\delta)(W)$$  we can suppose without loss
that $K<1.$

We distinguish then three cases corresponding to different signs
of $W-K \bar W.$  If  $W-K \bar W=0$ then one can take
$\varepsilon=1/2$ since
$$(\mu_1(\delta)-K\bar\mu_1(\delta))+
(\mu_2(\delta)-K\bar\mu_2(\delta))+
(\mu_3(\delta)-K\bar\mu_3(\delta))=0.$$

Let $W-K \bar W=W-\bar W+(1-K)\bar W<0.$ Then
$$(\mu_1(\delta)-K\bar\mu_1(\delta))+
(\mu_2(\delta)-K\bar\mu_2(\delta))+
(\mu_3(\delta)-K\bar\mu_3(\delta))=-3(W-K \bar W)\delta>0$$ and
$$
   \mu_3(\delta)-K\bar\mu_3(\delta)= \mu_3(\delta)-\bar\mu_3(\delta)
   +(1-K)\bar\mu_3(\delta)=\mu_3(\delta)(W')(W-\bar W)+(1-K)\bar\mu_3(\delta) $$
for $W'\in (W, \bar W).$ Since
$$\bar\mu_3(\delta)\le {\delta -3\over 3\sqrt 3} < {-1\over 3 \sqrt 3}
\le -\bar W,\;\;\mu'_3(\delta)(W')\le -2/3-\delta\le  -5/3 <-1$$
we get
$$\mu_3(\delta)-K\bar\mu_3(\delta)< -(W-\bar W+(1-K)\bar W)=
-(W-K\bar W)$$ and one can take
$\varepsilon=(2+3\delta)^{-1}=1/(2+3\delta)\ge
(2-\delta)/(4+\delta) $.

Let then $W-K \bar W=W-\bar W+(1-K)\bar W>0.$ We get
$$(\mu_1(\delta)-K\bar\mu_1(\delta))+
(\mu_2(\delta)-K\bar\mu_2(\delta))+
(\mu_3(\delta)-K\bar\mu_3(\delta))=-3(W-K \bar W)<0.$$ If $\bar
W\ge 0 $ then
$$\mu_2(\delta)-K\bar\mu_2(\delta)= \mu_2(\delta)-\bar\mu_2(\delta)
   +(1-K)\bar\mu_2(\delta)=\mu'_2(\delta)(W')(W-\bar W)+(1-K)
   \bar\mu_2(\delta)\ge $$ $$(2-\delta)(W-\bar W) +
   (1-K)(2-\delta)\bar W\ge (2-\delta)(W-K\bar W)$$
   which gives again $\varepsilon=(2-\delta)/(4+\delta).$

    Let $\bar W< 0,\; W\ge 0. $  Then
    $$\mu_2(\delta)-K\bar\mu_2(\delta)\ge (2-\delta) W +K(2-\delta)\bar
    W=(2-\delta)(W-K\bar W).$$
 Let finally $\bar W<  W< 0$.   Then the same inequality holds
 since the function $f(W):= \mu_2(\delta)(W)/W$ is decreasing for
 $W\in[{-1\over 3\sqrt3},0] $ and $f(0)=(2-\delta). $

 \medskip {\em End of
proof of Proposition 4.1.} Let us then   recall that
$$D^2w_{\delta }(a)_{|H}=(D^2P(a)-\delta P(a))_{|H}$$
for any    plane  $H$ orthogonal to a unit vector $a$. Applying
Corollary 4.1 to $H_{19}= a^{\perp}\bigcap b^{\perp}\bigcap H$
  and then  Lemma 4.1. with $K:= |b| ^{-\delta} $ we get the result in  all cases 
 except  $K=1, \; W(a)=W(b);$  but in this  exceptional case the trace of
 $H_{\delta}(a,b,O)$ vanishes and the claim is valid for $ \varepsilon={1\over 20}$.

\bigskip
 Proposition 4.1 and Lemma 2.1 give a proof of Theorem 1.1.
 Indeed,  we set $K$ to be the dual cone  $ K:=K_\lambda^{\ast}$ where
  $$  K_\lambda=
\{( \lambda_1,...  ,\lambda_{n})\in [C/\lambda, C\lambda]: {\hbox{
for some}}\; C>0 \;\}
$$ with $ n=21, \;  \lambda={1\over \varepsilon}.$
Then Proposition 4.1 gives  the $K-$cone condition in
 Lemma 2.1 on   $T_{\sigma_0 }\Lambda (B)$ for ${\sigma_0=id \in \Sigma_{21} }$
which implies the same  condition on the whole
 $M=\bigcup_{\sigma \in \Sigma_{21} }\   T_{\sigma }\Lambda (B) $ as well.

\section{Isaacs equation}

\medskip
We prove here Theorem 1.2.
Denote for $C>0$ by $K_C\subset S^2({\bf R}^2) $ the cone of positive symmetric matrix with
the ellipticity constant $C$, i.e., if $A\in K_C, \  A=\{ a_{ij} \}$ then
$$C^{-1}|\xi  |^2 \leq \sum a_{ij}\xi_i\xi_j \leq C|\xi |^2.$$

\medskip
{\bf Lemma 5.1.} {\it Let $C>0$ and let $w\in C^{\infty }({\bf R}^n\setminus 0) $
 be a homogeneous order $\alpha , 1<\alpha
\leq 2$ function. Assume that for any two points $x,y\in {\bf R}^n,\ 0<|x|,|y|\le 1$, there exists a matrix
$A\in K_C$ orthogonal to both forms $D^2w(x), D^2w(y),$ 
$$Tr(A D^2w(x))=Tr (AD^2w(y))=0.$$
Then $w$ is a viscosity solution to an Isaacs equation. }

\medskip
{\bf Proof.} Set
$$S=\{ a\in K_C,  tr \ a=1\}.$$
Denote
$$\Gamma =D^2w(S^{n-1}_1)\subset S^2({\bf R}^n).$$
Let
$$b\in S^2({\bf R}^n).$$
Denote
$$ B =\{ z\in S^2({\bf R}^n), zb>0 \},$$
$$b^*=B\cap S.$$

We define a two-parametric set of quadratic forms $L_{ab}\subset S^2({\bf R}^n)$
parametrized by $b\in \Gamma $ and $a\in b^*,\; a=\{ a_{ij}\}$. Denote by $L_{ab}$
the linear elliptic operator (1.5) with the coefficients $a_{ij}$ given by the parameter $a$.
 Then $L_{ab}$ is a uniformly elliptic operator with the ellipticity constant $C$. 
We are going to show that 
$$\sup_b \inf_a L_{ab} w =0.\leqno (5.1)$$
Let $x\in B,\ |x|\neq 0$. Choose $b=D^2w(x/|x|)$. Then since $D^2w(x)$ is proportional to $D^2w(x/|x|)$
we have
$$ \inf_{a\in z^*} L_{ab} w(x) =0.\leqno (5.2)$$
Assume now that $b_0\neq b$. By our assumptions there exists $A\in b_0^*\cap b^*$,
such that $Ab=Ab_0=0$. Thus
$$ \inf_{a\in b_0^*} L_{ab_{0}} z \leq 0.\leqno (5.3)$$
Now from (5.2) and (5.3) the equality (5.1) follows  immediately .

\medskip
 Recall that a symmetric matrix $A$ is  called strictly hyperbolic if
$$\frac{1}{ M} < -\frac{\lambda_1(A)}{\lambda_n (A)} < M$$
for a positive $M$.
To finish the   proof we note that the results of Section 4 imply that the 
form  $ \alpha F_1{D^2w_{\delta}} _{|H}(x) -\beta  F_2{D^2w_{\delta}} _{|H}(y) $
 is strictly hyperbolic for positive $ \alpha, \beta$; since the function $w_{\delta}$ is 
odd, it remains true for any  $ (\alpha, \beta) \in {\bf R}^2\setminus \{0\}$ and we can apply
 the following result.  

\medskip
{\bf Lemma 5.2.} {\it Let $F_1, F_2$ be two quadratic forms  in ${\bf R}^n$  
s.t. the form \linebreak $ \alpha F_1+\beta  F_2 $
 is strictly  hyperbolic for any  $ (\alpha, \beta) \in {\bf R}^2\setminus \{0\}$. 
Then there exist $C>0$ and a positive  quadratic form $Q\in K_C $ orthogonal 
to both forms $F_1, F_2$, 
$$Tr(F_1Q)=Tr(F_2Q)=0. $$ }
\medskip
{\em Proof.} We  can suppose w.r.g. that $F_1$ is traceless, $Tr(F_1)=0. $ Let $D\subset S_1^{n-1}$ be a 
minimal (with respect to inclusion) domain on which $F_1$ does not change the sign. We can assume  w.r.g.
that $F_1\mid_D>0.$  Our first claim is that $F_2$ changes the sign on the border $\partial D$ of $D$. Indeed,
 if not we assume w.r.g. $F_2\mid_{\partial D}\ge 0. $ Let for $ t\in [0,1]$ define $ D_t$ as the union of the connected 
components of the set $\{x\in S_1^{n-1}: \cos(\pi t) F_1(x)+\sin(\pi t)  F_2(x)>0 \}$ with non-void intersection with $D$,
thus  $  D_0=D$. If for some $ s\in [0,1]$ we get $ D_s \bigcap D\neq D_s$ we are  done and we can thus assume that
$ \forall s\in [0,1],\; D_s\subset D.$ If for some $ s\in [0,1]$ the set  $ D_s$ becomes empty, then there is $s'\in [0,s[$ s.t.
$\bar D_s'$ contains an isolated point  $x_0$ with $\cos(\pi s') F_1(x_0)+\sin(\pi s ')  F_2(x_0)=0$ which is impossible
 since then 0 would be a maximal eigenvalue of the strictly hyperbolic form  $\cos(\pi s') F_1+\sin(\pi s ')  F_2.$ In particular,
$D_1$ is non-empty which is impossible since $F_1=-F_0.$ 

Since $F_2$ changes the sign on the border $\partial D$  of $D$, there exist two points $a_1, a_2\in \partial D$ 
with $F_1(a_1)= F_1(a_2)=0, F_2(a_1)= a>0, F_2(a_2)=-a.$  Let $m= Tr(F_2 (\sum x^2_i)),$ changing 
the sign if necessary we can suppose that  $m\ge 0$. If $m=0$ we are done with $Q=\sum x^2_i$, 
thus we suppose $m>0$.
Then the form $Q_0(x):=(x,a_2)^2$ is clearly 
orthogonal to $F_1$ and one has $Tr(F_2 Q_{0})=-a.$ Let $l:=a/m>0$. Then the form 
$ Q_{0,l}:=Q_0(x)+l\sum x^2_i$ is positive, orthogonal to $F_1$ since $F_1$  is traceless,
 and one has  $Tr(F_2 Q_{0,l})=-a+ml=0.$
One notes then that the ellipticity constant of the form $ Q_{0,l}$ depends (upper semi-) continuously on
 $ (F_1, F_2)$, thus its maximum $C$ on $S ^{n-1}_1\times S ^{n-1}_1$ is finite.
 
  The lemma is proved.

\bigskip
 \centerline{REFERENCES}

\medskip
\noindent [Ad] J.F. Adams, {\it Lectures on Exceptional Lie Groups,} Univ. Chicago Press, Chicago, 1996.

\medskip
\noindent [Bae] J. Baez, {\it Octonions}, Bull. Amer. Math. Soc.,
39 (2002), 145--205.

\medskip
 \noindent [C] L. Caffarelli,  {\it Interior a priory estimates for solutions
 of fully nonlinear equations}, Ann. Math. 130 (1989), 189--213.

\medskip
 \noindent [CC] L. Caffarelli, X. Cabre, {\it Fully Nonlinear Elliptic
Equations}, Amer. Math. Soc., Providence, R.I., 1995.

\medskip
 \noindent [CIL]  M.G. Crandall, H. Ishii, P-L. Lions, {\it User's
guide to viscosity solutions of second order partial differential
equations,} Bull. Amer. Math. Soc. (N.S.), 27(1) (1992), 1--67.

\medskip
 \noindent [CNS] L. Caffarelli, L. Nirenberg, J. Spruck, {\it The Dirichlet
 problem for nonlinear second order elliptic equations III. Functions
  of the eigenvalues of the Hessian, } Acta Math.
   155 (1985), no. 3-4, 261--301.

   \medskip
 \noindent [F] A. Friedman, {\it Differential games }, Pure and Applied Mathematics,
 vol. 25, John Wiley and Sons, New York, 1971

\medskip
 \noindent [GT] D. Gilbarg, N. Trudinger, {\it Elliptic Partial
Differential Equations of Second Order, 2nd ed.}, Springer-Verlag,
Berlin-Heidelberg-New York-Tokyo, 1983.

\medskip
\noindent [NV1] N. Nadirashvili, S. Vl\u adu\c t, {\it
Nonclassical solutions of fully nonlinear elliptic equations,}
Geom. Func. An. 17 (2007), 1283--1296.

\medskip
\noindent [NV2] N. Nadirashvili, S. Vl\u adu\c t, {\it Singular Viscosity Solutions to Fully Nonlinear 
Elliptic Equations},  J. Math.Pures Appl., 89 (2008), 107-113.

\medskip
\noindent [T1] N. Trudinger, {\it Weak solutions of Hessian
equations,} Comm. Partial Differential Equations 22 (1997), no.
7-8, 1251--1261.

\medskip
\noindent [T2] N. Trudinger, {\it On the Dirichlet problem for
Hessian equations,} Acta Math.
   175 (1995), no. 2, 151--164.

\medskip
\noindent [T3] N. Trudinger, {\it H\"older gradient estimates for
fully nonlinear elliptic equations,} Proc. Roy. Soc. Edinburgh
Sect. A 108 (1988), 57--65.

\medskip
\noindent [We] G. Weyl, {\it Das asymptotische Verteilungsgezets
des Eigenwerte lineare partieller Differentialgleichungen,} Math.
Ann.
   71 (1912), no. 2,  441--479.
\end{document}